\documentclass[a4paper,reqno,12pt]{amsart}
\usepackage{amsmath,amssymb,amsbsy,amsfonts,amsthm,latexsym,
            amsopn,amstext,amsxtra,euscript,amscd, color}
\usepackage[colorinlistoftodos]{todonotes}
\usepackage[colorlinks=true]{hyperref}
\hypersetup{urlcolor=blue, citecolor=green}

\usepackage{color,transparent}
\usepackage{hyperref}   
\usepackage[margin=3cm]{geometry}

\begin{document}
\bibliographystyle{plain}
\newfont{\teneufm}{eufm10}
\newfont{\seveneufm}{eufm7}
\newfont{\fiveeufm}{eufm5}
%
%
\newfam\eufmfam
              \textfont\eufmfam=\teneufm \scriptfont\eufmfam=\seveneufm
              \scriptscriptfont\eufmfam=\fiveeufm
\def\bbbr{{\rm I\!R}}
\def\bbbm{{\rm I\!M}}
\def\bbbn{{\rm I\!N}}
\def\bbbf{{\rm I\!F}}
\def\bbbh{{\rm I\!H}}
\def\bbbk{{\rm I\!K}}
\def\bbbp{{\rm I\!P}}
\def\bbbone{{\mathchoice {\rm 1\mskip-4mu l} {\rm 1\mskip-4mu l}
{\rm 1\mskip-4.5mu l} {\rm 1\mskip-5mu l}}}
\def\bbbc{{\mathchoice {\setbox0=\hbox{$\displaystyle\rm C$}\hbox{\hbox
to0pt{\kern0.4\wd0\vrule height0.9\ht0\hss}\box0}}
{\setbox0=\hbox{$\textstyle\rm C$}\hbox{\hbox
to0pt{\kern0.4\wd0\vrule height0.9\ht0\hss}\box0}}
{\setbox0=\hbox{$\scriptstyle\rm C$}\hbox{\hbox
to0pt{\kern0.4\wd0\vrule height0.9\ht0\hss}\box0}}
{\setbox0=\hbox{$\scriptscriptstyle\rm C$}\hbox{\hbox
to0pt{\kern0.4\wd0\vrule height0.9\ht0\hss}\box0}}}}
\def\bbbq{{\mathchoice {\setbox0=\hbox{$\displaystyle\rm
Q$}\hbox{\raise
0.15\ht0\hbox to0pt{\kern0.4\wd0\vrule height0.8\ht0\hss}\box0}}
{\setbox0=\hbox{$\textstyle\rm Q$}\hbox{\raise
0.15\ht0\hbox to0pt{\kern0.4\wd0\vrule height0.8\ht0\hss}\box0}}
{\setbox0=\hbox{$\scriptstyle\rm Q$}\hbox{\raise
0.15\ht0\hbox to0pt{\kern0.4\wd0\vrule height0.7\ht0\hss}\box0}}
{\setbox0=\hbox{$\scriptscriptstyle\rm Q$}\hbox{\raise
0.15\ht0\hbox to0pt{\kern0.4\wd0\vrule height0.7\ht0\hss}\box0}}}}
\def\bbbt{{\mathchoice {\setbox0=\hbox{$\displaystyle\rm.
T$}\hbox{\hbox to0pt{\kern0.3\wd0\vrule height0.9\ht0\hss}\box0}}
{\setbox0=\hbox{$\textstyle\rm T$}\hbox{\hbox
to0pt{\kern0.3\wd0\vrule height0.9\ht0\hss}\box0}}
{\setbox0=\hbox{$\scriptstyle\rm T$}\hbox{\hbox
to0pt{\kern0.3\wd0\vrule height0.9\ht0\hss}\box0}}
{\setbox0=\hbox{$\scriptscriptstyle\rm T$}\hbox{\hbox
to0pt{\kern0.3\wd0\vrule height0.9\ht0\hss}\box0}}}}
\def\bbbs{{\mathchoice
{\setbox0=\hbox{$\displaystyle     \rm S$}\hbox{\raise0.5\ht0\hbox
to0pt{\kern0.35\wd0\vrule height0.45\ht0\hss}\hbox
to0pt{\kern0.55\wd0\vrule height0.5\ht0\hss}\box0}}
{\setbox0=\hbox{$\textstyle        \rm S$}\hbox{\raise0.5\ht0\hbox
to0pt{\kern0.35\wd0\vrule height0.45\ht0\hss}\hbox
to0pt{\kern0.55\wd0\vrule height0.5\ht0\hss}\box0}}
{\setbox0=\hbox{$\scriptstyle      \rm S$}\hbox{\raise0.5\ht0\hbox
to0pt{\kern0.35\wd0\vrule height0.45\ht0\hss}\raise0.05\ht0\hbox
to0pt{\kern0.5\wd0\vrule height0.45\ht0\hss}\box0}}
{\setbox0=\hbox{$\scriptscriptstyle\rm S$}\hbox{\raise0.5\ht0\hbox
to0pt{\kern0.4\wd0\vrule height0.45\ht0\hss}\raise0.05\ht0\hbox
to0pt{\kern0.55\wd0\vrule height0.45\ht0\hss}\box0}}}}
\def\bbbz{{\mathchoice {\hbox{$\sf\textstyle Z\kern-0.4em Z$}}
{\hbox{$\sf\textstyle Z\kern-0.4em Z$}}
{\hbox{$\sf\scriptstyle Z\kern-0.3em Z$}}
{\hbox{$\sf\scriptscriptstyle Z\kern-0.2em Z$}}}}
\def\ts{\thinspace}

\newtheorem{theorem}{Theorem}
\newtheorem{lem}{Lemma}
\newtheorem{lemma}[lem]{Lemma}

\newtheorem{claim}[theorem]{Claim}
\newtheorem{cor}[theorem]{Corollary}
\newtheorem{prop}[theorem]{Proposition}
\newtheorem{definition}[theorem]{Definition}
\newtheorem{remark}[theorem]{Remark}
\newtheorem{question}[theorem]{Open Question}
\newtheorem{example}[theorem]{Example}
\newtheorem{problem}[theorem]{Problem}

\def\qed{\ifmmode
\squareforqed\else{\unskip\nobreak\hfil
\penalty50\hskip1em\null\nobreak\hfil\squareforqed
\parfillskip=0pt\finalhyphendemerits=0\endgraf}\fi}

\def\squareforqed{\hbox{\rlap{$\sqcap$}$\sqcup$}}

\def \C {{\mathbb C}}
\def \F {{\mathbb F}}
\def \L {{\mathbb L}}
\def \K {{\mathbb K}}
\def \Q {{\mathbb Q}}
\def \Z {{\mathbb Z}}
\def\cA{{\mathcal A}}
\def\cB{{\mathcal B}}
\def\cC{{\mathcal C}}
\def\cD{{\mathcal D}}
\def\cE{{\mathcal E}}
\def\cF{{\mathcal F}}
\def\cG{{\mathcal G}}
\def\cH{{\mathcal H}}
\def\cI{{\mathcal I}}
\def\cJ{{\mathcal J}}
\def\cK{{\mathcal K}}
\def\cL{{\mathcal L}}
\def\cM{{\mathcal M}}
\def\cN{{\mathcal N}}
\def\cO{{\mathcal O}}
\def\cP{{\mathcal P}}
\def\cQ{{\mathcal Q}}
\def\cR{{\mathcal R}}
\def\cS{{\mathcal S}}
\def\cT{{\mathcal T}}
\def\cU{{\mathcal U}}
\def\cV{{\mathcal V}}
\def\cW{{\mathcal W}}
\def\cX{{\mathcal X}}
\def\cY{{\mathcal Y}}
\def\cZ{{\mathcal Z}}
\newcommand{\rmod}[1]{\: \mbox{mod}\: #1}

\def\tcN{\cN^\mathbf{c}}
\def\F{\mathbb F}
\def\Tr{\operatorname{Tr}}
\def\mand{\qquad \mbox{and} \qquad}
\renewcommand{\vec}[1]{\mathbf{#1}}
\def\eqref#1{(\ref{#1})}
\newcommand{\ignore}[1]{}
\hyphenation{re-pub-lished}
\parskip 1.5 mm
\def\lln{{\mathrm Lnln}}
\def\Res{\mathrm{Res}\,}
\def\F{{\bbbf}}
\def\Fp{\F_p}
\def\fp{\Fp^*}
\def\Fq{\F_q}
\def\ff{\F_2}
\def\ffn{\F_{2^n}}
\def\K{{\bbbk}}
\def \Z{{\bbbz}}
\def \N{{\bbbn}}
\def\Q{{\bbbq}}
\def \R{{\bbbr}}
\def \P{{\bbbp}}
\def\Zm{\Z_m}
\def \Um{{\mathcal U}_m}
\def \Bf{\frak B}
\def\Km{\cK_\mu}
\def\va {{\mathbf a}}
\def \vb {{\mathbf b}}
\def \vc {{\mathbf c}}
\def\vx{{\mathbf x}}
\def \vr {{\mathbf r}}
\def \vv {{\mathbf v}}
\def\vu{{\mathbf u}}
\def \vw{{\mathbf w}}
\def \vz {{\mathbfz}}
\def\\{\cr}
\def\({\left(}
\def\){\right)}
\def\fl#1{\left\lfloor#1\right\rfloor}
\def\rf#1{\left\lceil#1\right\rceil}
\def\flq#1{{\left\lfloor#1\right\rfloor}_q}
\def\flp#1{{\left\lfloor#1\right\rfloor}_p}
\def\flm#1{{\left\lfloor#1\right\rfloor}_m}
\def\Al{{\sl Alice}}
\def\Bob{{\sl Bob}}
\def\Or{{\mathcal O}}
\def\inv#1{\mbox{\rm{inv}}\,#1}
\def\invM#1{\mbox{\rm{inv}}_M\,#1}
\def\invp#1{\mbox{\rm{inv}}_p\,#1}
\def\Ln#1{\mbox{\rm{Ln}}\,#1}
\def \nd {\,|\hspace{-1.2mm}/\,}
\def\ord{\mu}
\def\E{\mathbf{E}}
\def\Cl{{\mathrm {Cl}}}
\def\epp{\mbox{\bf{e}}_{p-1}}
\def\ep{\mbox{\bf{e}}_p}
\def\eq{\mbox{\bf{e}}_q}
\def\bm{\bf{m}}
\newcommand{\floor}[1]{\lfloor {#1} \rfloor}
\newcommand{\comm}[1]{\marginpar{
\vskip-\baselineskip
\raggedright\footnotesize
\itshape\hrule\smallskip#1\par\smallskip\hrule}}
\def\rem{{\mathrm{\,rem\,}}}
\def\dist {{\mathrm{\,dist\,}}}
\def\etal{{\it et al.}}
\def\ie{{\it i.e. }}
\def\veps{{\varepsilon}}
\def\eps{{\eta}}
\def\ind#1{{\mathrm {ind}}\,#1}
               \def \MSB{{\mathrm{MSB}}}
\newcommand{\abs}[1]{\left| #1 \right|}


\author[B. EDJEOU and B. FAYE]{Bilizimb\'ey\'e EDJEOU and Bernadette Faye}



\address{Bilizimb\'ey\'e EDJEOU\newline
         \indent Ecole Sup\'erieure d'Informatique et de Gestion \newline
          \indent  ESIG GLOBAL SUCCESS \newline
        \indent Lom\'e 149\newline
         \indent Togo}
\email{olivieredjeou@gmail.com}

\address{Bernadette Faye\newline
         \indent UFR SATIC, Universit\'e Alioune Diop de Bambey \newline
         \indent Departement de Mathematiques\newline
       \indent Bambey 30\newline
         \indent Diourbel, S\'en\'egal}
         
\email{bernadette.faye@uadb.edu.sn}



\title{Pell and Pell-Lucas Numbers as difference of two repdigits}

\maketitle
\begin{abstract}
Let $ \{P_{n}\}_{n\geq 0} $ be the sequence of Pell numbers  defined by $ P_0=0 $, $ P_1 =1$ and $ P_{n+2}= 2P_{n+1} +P_n$ for all $ n\geq 0 $ and let $ \{Q_{n}\}_{n\geq 0} $ be its companion sequence, the Pell-Lucas numbers  defined by $ Q_0=Q_1 =2$ and $ Q_{n+2}= 2Q_{n+1} +Q_n$ for all $ n\geq 0 $ . In this paper, we find all Pell and Pell-Lucas numbers which can be written as difference of two repdigits. It is shown that the largest Pell and Pell-Lucas numbers which can be written as difference of two repdigits are 

$$P_6=70= 77-7 \quad\quad \hbox{and} \quad\quad Q_7 = 478=555-77.$$
\end{abstract}

\maketitle
\section{Introduction}

Let $(P_n)_{n\ge 0}$ be the sequence of Pell numbers given by $P_0=0$, $P_1=1$ and 
$$P_{n+2}= 2P_{n+1} + P_n \quad {\text{\rm for all}}\quad n\geq 0.
$$
The Pell-Lucas sequence $(Q_n)_{n\ge 0}$ satisfies the same recurrence as the sequence of Pell numbers 
with initial conditions $Q_0=Q_1=2.$ If $(\alpha,\beta)=(1+{\sqrt{2}},1-{\sqrt{2}})$ is the pair of roots of the characteristic equation $x^2-2x-1=0$ 
of both the Pell and Pell-Lucas numbers, then the Binet formulas for their general terms are:
$$
P_n=\frac{\alpha^n-\beta^n}{\alpha-\beta}\quad {\text{\rm and}}\quad Q_n=\alpha^n+\beta^n\quad {\text{\rm for~all}}\quad n\ge 0.
$$ 
This implies easily that the inequalities
\begin{equation}
\label{eq:sizePn}
\alpha^{n-2}\le P_n\le  \alpha^{n-1}\quad\hbox{for $n\geq 2$}
\end{equation}
and 
\begin{equation}
\label{eq:siezQn}
\alpha^{n-1}\le Q_n\le \alpha^{n+1}\quad\hbox{for $n\geq 1$}
\end{equation}
hold.

Given an integer $g>1$, a base $g$-repdigit is a number of the form 

$$N=a\cdot\frac{g^m-1}{g-1} \quad\quad \hbox{for some $m\geq 1$ and $a\in\{1,\ldots,(g-1)\}$}.$$

When $g=10$, such number are better know as  a \textit{repdigit}.
Investigation of the repdigits in the second-order linear recurrence sequences has been of interest to mathematicians. All Pell  and Pell-Lucas numbers which are repdigits have been found in \cite{FL1}. The largest repdigits in Pell and Pell-Lucas sequences are $P_3 = 5$ and $Q_2 = 6$, respectively. Subsequently,  Fibonacci, Lucas, Pell and Pell-Lucas numbers which are expressible as sum of two repdigits have been studied in [\cite{ALT1},\cite{ALT2},\cite{AL}].

Then, it is natural to replace in the above problems the sums by diferences. Recently, Erduwan at al. found in \cite{EKL} all Fibonacci and Lucas numbers which are difference of two repdigits. 

In this paper, we study the same question and ask which
Pell and Pell-Lucas numbers can be written as a difference of two repdigits.  That is, we study the Diophantine equations

\noindent
\begin{eqnarray}\label{Main1a}
P_k=a_1\cdot\frac{10^n-1}{9}-a_2\cdot\frac{10^m-1}{9}  
\end{eqnarray}
and
\begin{eqnarray}\label{Main1b}
Q_k=a_1\cdot\frac{10^n-1}{9} - a_2\cdot\frac{10^m-1}{9}  
\end{eqnarray}
where $k$, $n$, $m$ are some positive integers with $n\geq 2$ and  $a_1,a_2\in\{1,\ldots,9\}$.
We have the following results.


\begin{theorem}\label{Main}
If $P_k$ is expressible as a difference of two repdigits, then

$$P_k\in\{2, 5, 12, 29, 70\},$$
with

\begin{eqnarray*}
P_2 &=&2=11-9, \quad P_3=5=11-6,\\
P_4 &=& 12=111-99, \quad P_5=29=33-4,\\
P_6 &=& 70=77-7.
\end{eqnarray*}
\end{theorem}

\begin{theorem}\label{Main1}
If $Q_k$ is expressible as a difference of two repdigits, then

$$Q_k\in\{2, 5, 6, 14, 34, 82, 478\},$$
with

\begin{eqnarray*}
Q_0 &=&2=11-9, \quad Q_1=5=11-6,\\
Q_2 &=& 6=11-5, \quad Q_3=14=22-8,\\
Q_4 &= &34=111-77, \quad Q_5=82=88-6,\\
Q_7 &=& 478=555-77. 
\end{eqnarray*}
\end{theorem}

\section{Preliminary results}
\subsection{Linear forms in logarithms}
To prove our main result Theorem \ref{Main}, we use several times a Baker--type lower bound for a nonzero linear form in logarithms of algebraic numbers. There are many such  bounds in the literature  like that of Baker and W{\"u}stholz from \cite{bawu07}. In this paper we use the result of Matveev \cite{MatveevII}, which is one of our main tools. We start
with recalling some basic definitions and results from algebraic number theory

\noindent
Let $ \gamma $ be an algebraic number of degree $ d $ with minimal primitive polynomial over the integers
$$ a_{0}x^{d}+ a_{1}x^{d-1}+\cdots+a_{d} = a_{0}\prod_{i=1}^{d}(x-\gamma^{(i)}),$$
where the leading coefficient $ a_{0} $ is positive and the $ \eta^{(i)} $'s are the conjugates of $ \gamma $. Then the \textit{logarithmic height} of $ \gamma $ is given by
$$ h(\gamma) := \dfrac{1}{d}\left( \log a_{0} + \sum_{i=1}^{d}\log\left(\max\{|\gamma^{(i)}|, 1\}\right)\right).$$

In particular, if $ \gamma = p/q $ is a rational number with $ \gcd (p,q) = 1 $ and $ q>0 $, then $ h(\gamma) = \log\max\{|p|, q\} $. The following are some of the properties of the logarithmic height function $ h(\cdot) $, which will be used in the next sections of this paper without reference:
\begin{eqnarray}
\label{eq:height}
h(\eta\pm \gamma) &\leq& h(\eta) +h(\gamma) +\log 2,\nonumber\\
h(\eta\gamma^{\pm 1})&\leq & h(\eta) + h(\gamma),\\
h(\eta^{s}) &=& |s|h(\eta) ~~~~~~ (s\in\mathbb{Z}). \nonumber
\end{eqnarray}

\begin{theorem}[Matveev]\label{Matveev11} Let $\gamma_1,\ldots,\gamma_t$ be positive real algebraic numbers in a real algebraic number field 
$\mathbb{K}$ of degree $D$, $b_1,\ldots,b_t$ be nonzero integers, and assume that
\begin{equation}
\label{eq:Lambda}
\Lambda:=\gamma_1^{b_1}\cdots\gamma_t^{b_t} - 1,
\end{equation}
is nonzero. Then
$$
\log |\Lambda| > -1.4\times 30^{t+3}\times t^{4.5}\times D^{2}(1+\log D)(1+\log B)A_1\cdots A_t,
$$
where
$$
B\geq\max\{|b_1|, \ldots, |b_t|\},
$$
and
$$A
_i \geq \max\{Dh(\gamma_i), |\log\gamma_i|, 0.16\},\qquad {\text{for all}}\qquad i=1,\ldots,t.
$$
\end{theorem} 
\subsection{Baker-Davenport reduction lemma}
During the calculations, we get upper bounds on our variables which are too large, thus we need to reduce them. To do so, we use some results from the theory of continued fractions. Specifically, for a nonhomogeneous linear form in two integer variables, we use a slight variation of a result due to Dujella and Peth{\H o} (see \cite{dujella98}, Lemma 5a), which  is itself a generalization of a result of Baker and Davenport \cite{BD69}. For a real number $X$, we write  $||X||:= \min\{|X-n|: n\in\mathbb{Z}\}$ for the distance from $X$ to the nearest integer.
\begin{lemma}[Dujella, Peth\H o]\label{Dujjella}
Let $M$ be a positive integer, $p/q$ be a convergent of the continued fraction of the irrational number $\tau$ such that $q>6M$, and  $A,B,\mu$ be some real numbers with $A>0$ and $B>1$. Let further 
$\varepsilon: = ||\mu q||-M||\tau q||$. If $ \varepsilon > 0 $, then there is no solution to the inequality
$$
0<|u\tau-v+\mu|<AB^{-w},
$$
in positive integers $u,v$ and $w$ with
$$ 
u\le M \quad {\text{and}}\quad w\ge \dfrac{\log(Aq/\varepsilon)}{\log B}.
$$
\end{lemma}


\noindent
Finally, the following lemma is also useful. It it can be found in \cite{weber}.

\begin{lemma}
\label{weber}
If $a,x\in \mathbb{R}$. If $0<a<1$ and $|x|<a.$ Then, 
$$|\log(x+1)|< \frac{-\log(1-a)}{a}\cdot |x|$$
and 

$$|x|<\frac{a}{1-e^{-a}}\cdot |e^x-1|.$$
\end{lemma}

\section{Proof of Theorem \ref{Main}}
\noindent
Assume that the equation  \eqref{Main1a} holds. Let $1\leq k \leq 149$ and $n\geq 2$. Then, by using Sagemath, we obtain only the solutions listed in Theorem \ref{Main}.
\medskip

From now, we assume that $k\geq 150.$ If $n=m$, then it follows that $a_1>a_2$, which means that $P_k$ is a repdigit. But the largest repdigit in $P_k$ is $5$. Thus we get a contradiction since $k\geq 150.$ 
\medskip


So, we assume now that  $k\geq 150$, $n-m\geq 1$. Therefore, using \eqref{eq:sizePn}, we obtain the  inequality

$$\frac{\alpha^{2n}}{20}<\frac{10^{n-1}}{2}<10^{n-1}-10^{m}<a_1\frac{10^n-1}{9}-a_2\frac{10^m-1}{9}=P_k <\alpha^{k-1}, $$
which implies that $2n<k+3.$ In particular $n \leq k.$ On the other hand, we rewrite equation \eqref{Main1a} as

$$\frac{\alpha^k - \beta^k}{2\sqrt{2}}= a_1\frac{10^n-1}{9}-a_2\frac{10^m-1}{9}$$
to obtain

\begin{equation}
\label{eq:1}
\frac{9\alpha^k}{2\sqrt{2}}-a_110^n=\frac{9\beta^k}{2\sqrt{2}}-a_210^m-(a_1-a_2).
\end{equation}
Taking absolute value of both sides of equation \eqref{eq:1}, we obtain 
\begin{equation}
\label{eq:2}
\Big|\frac{9\alpha^k}{2\sqrt{2}}-a_110^n\Big|\leq \frac{9|\beta|^k}{2\sqrt{2}}+a_210^m+|a_1-a_2|.
\end{equation}
Divinding both sides of \eqref{eq:2} by $a_110^n$, we obtain

\begin{align*}
\Big|\frac{9\cdot 10^{-n}\alpha^k}{a_12\sqrt{2}}-1\Big|& \leq \frac{9|\beta|^k}{a_110^n2\sqrt{2}}+\frac{a_210^m}{a_110^n}+\frac{|a_1-a_2|}{a_110^n}\\
&\leq \frac{9|\beta|^k}{10^{n-m+1}2\sqrt{2}}+\frac{9}{10^{n-m}}+\frac{8}{10^{n-m+1}}.
\end{align*}
This implies that 

\begin{equation}
\label{eq:3}
\Big|\frac{9\cdot 10^{-n}\alpha^k}{a_12\sqrt{2}}-1\Big| \leq \frac{9.81}{10^{n-m}}\leq \frac{9.81}{10}=0.981,
\end{equation}
for $n-m\geq 1.$
Now we apply Theorem \ref{Matveev11} to the left-hand side of the above inequality with $(\gamma_1,\gamma_2,\gamma_3)=(\alpha,10,9/a_12\sqrt{2})$ and $(b_1,b_2,b_3)=(k,-n,1)$. Note that $\gamma_1,\gamma_2$ and $\gamma_3$ are positive real numbers and elements of the field $\mathbb{K}=\mathbb{Q}(\sqrt{2})$. Therefore the degree of the field $\mathbb{K}$ is equal to $D=2$.
Put 

$$\Lambda_1=\Big|\frac{9\cdot 10^{-n}\alpha^k}{a_12\sqrt{2}}-1\Big\vert.$$
If $\Lambda_1=0$, then $\alpha^k=10^na_12\sqrt{2}/9$. Conjugating in $\mathbb{K}$ gives that $\beta^k=-10^na_12\sqrt{2}/9$ and so $Q_k=\alpha^k+\beta^k=0$ which is a contradiction.  Therefore $\Lambda_1\neq 0$. Using the properties of the logarithmic height in \eqref{eq:height}, we have that 

$$h(\gamma_1)=h(\alpha)=\frac{\log\alpha}{2}\leq\frac{0.882}{2}, \quad h(\gamma_2)=\log10$$
and 

$$h(\gamma_3)=h(9/a_12\sqrt{2})< h(9)+h(a_1)+h(2\sqrt{2})<5.44$$ 
with $h(2\sqrt{2})=\log8/2$. We can take $A_1=\log(\alpha)=0.882, A_2=4.7$ and $A_3=10.9$ Since $n<k$ and $B\geq \max(n,k,1)$ we can take $B=k$. Thus, we obtain by inequality \eqref{eq:3} and  Theorem \ref{Matveev11} that 

$$(10.12)\cdot 10^{m-n}>\Lambda_1>exp(C\cdot (1+\log2)\cdot (1+\log(k))(0.882)(4.6)(10.9)),$$
where $C=-1.4\cdot 30^6\cdot 3^{4.5}\cdot 2^2.$ It follows that 

\begin{equation}
\label{eq:4}
(m-n)\log10< 4.29\cdot 10^{13}\cdot (1+\log(k))+\log(9.81).
\end{equation}

Rearranging equation \eqref{Main1a} 	as 

\begin{equation}
\label{eq:5}
\frac{\alpha^k}{2\sqrt{2}}-\frac{a_110^n-a_210^m}{9}=\frac{\beta^k}{2\sqrt{2}}-\frac{(a_1-a_2)}{9}
\end{equation}
and taking absolute value of both sides of \eqref{eq:5}, we get

\begin{equation}
\label{eq:6}
\Big|\frac{\alpha^k}{2\sqrt{2}}-\frac{a_110^n-a_210^m}{9}\Big|\leq \frac{\Big|\beta\Big|^k}{2\sqrt{2}}+\frac{\Big|a_1-a_2\Big|}{9}.
\end{equation}
Dividing both sides of the above inequality by $\frac{\alpha^k}{2\sqrt{2}}$ we obtain 

\begin{equation}
\label{eq:7}
\Big|1-\frac{(a_1-a_210^{m-n})\cdot 10^n\cdot \alpha^{-k}}{9\cdot 2\sqrt{2}}\Big|\leq \frac{1}{\alpha^{2k}}+\frac{16\sqrt{2}}{9\alpha^k}<\frac{4}{\alpha^k}
\end{equation}
Now, we can apply again Theorem \ref{Matveev11} to the above inequality  with 
$$(\gamma_1,\gamma_2,\gamma_3)=(\alpha,10,(a_1-a_210^{m-n})/18\sqrt{2})$$ and 

$$(b_1,b_2,b_3)=(-k,n,1)$$.

 Note that $\gamma_1,\gamma_2$ and $\gamma_3$ are positive real numbers and elements of the field $\mathbb{K}=\mathbb{Q}(\sqrt{2})$. Therefore the degree of the field $\mathbb{K}$ is equal to $D=2$. Let 
 
 $$\Lambda_2= 1-\frac{(a_1-a_210^{m-n})\cdot 10^n\cdot \alpha^{-k}}{18\sqrt{2}}.$$
If $\Lambda=0$, then get that $\alpha^{2k}\in \mathbb{Q}$ which is false for $k>0$. Using the properties of the logarithmic height in \eqref{eq:height}, we have that 

$$h(\gamma_1)=h(\alpha)=\frac{\log\alpha}{2}\leq\frac{0.882}{2}, \quad h(\gamma_2)=\log10$$
and 

\begin{align*}
h(\gamma_3)&=h((a_1-a_210^{m-n})/18\sqrt{2})\\
&\leq h(18)+h(\sqrt{2})+h(a_1)+h(a_2)+(n-m)h(10)+\log2\\
&\leq 8.33 +(n-m)\log10.
\end{align*}
So we can take $A_1=\log(\alpha)=0.9, A_2=4.7$ and $A_3=16.66+2(n-m)\log10$. Since $n<k$ and $B\geq \max(n,k,1)$ we can take $B=k$. Thus, we obtain by inequality \eqref{eq:7} and  Theorem \ref{Matveev11} that 

$$4\cdot \alpha^{-k}>\Lambda_2>\exp(C\cdot (1+\log2)\cdot (1+\log(k))(0.9)(4.6)(16.66+2(n-m)\log10)),$$
where $C=-1.4\cdot 30^6\cdot 3^{4.5}\cdot 2^2.$ By a simple computation, it follows that 

\begin{equation}
\label{eq:8}
k\log\alpha-\log 4< 4.1\cdot 10^{12}\cdot (1+\log(k))(16.66+2(n-m)\log 10).
\end{equation}
Now using inequalities \eqref{eq:4} and \eqref{eq:8}, a computer search with SageMath gives us that $k<1.25\cdot 10^{29}.$

\subsection{Reducing the Bounding of $ k $}
Now let us reduce the upper bound on $n$ by using the Baker-Davenport algorithm given in Lemma \ref{Dujjella}. Let 

$$z_1:= k\log \alpha-n\log 10 +\log(9/a_12\sqrt{2}).$$
From the inequality \eqref{eq:3}, we have that 
$$|z_1|=|e^{z_1}-1|\leq \frac{9.81}{10^{n-m}}\leq \frac{9.81}{10} = 0.981$$
for $n-m \geq 1.$ Choosing $a:=0.1$, we get the inequality

$$z_1=|\log(x+1)|\leq\frac{-\log(1-0.981)}{0.981}\cdot\frac{9.81}{10^{n-m}}<(10.67)\cdot 10^{m-n}$$
by Lemma \ref{weber}. Thus it follows that 

$$0<|k\log \alpha-n\log 10 +\log(9/a_12\sqrt{2})|<10.67\cdot 10^{m-n}.$$
\medskip

Dividing this inequality by $\log 10$, we obtain

\begin{equation}
\label{eq:9}
\Big|k\left(\frac{\log \alpha}{\log 10}\right)-n +\left(\frac{\log(9/a_12\sqrt{2})}{\log 10}\right)\Big|<4.7\cdot 10^{m-n}.
\end{equation}
We can take $\gamma:=\frac{\log \alpha}{\log 10}\not\in \mathbb{Q}$ and $M= 1.25\cdot 10^{29}.$ Then we found that $q_{67}$, the denominator of the $67$-th convergent of $\gamma$ is greater than $6M$. Now take 

$$\mu:=\frac{\log(9/a_12\sqrt{2})}{\log 10}.$$
In this case, considering the fact that $1\leq a_1\leq 9$, a quick computation with Sagemath gives us that 

$$\epsilon(\mu):=||\mu q_{68}||-M||\gamma q_{68}||>0.1.$$
Let $A=4.7$, $B=10$ and $w=n-m$ in Lemma \ref{Dujjella}. Thus, using Sagemath, we can say that the inequality \eqref{eq:9} has no solution if 

$$ \dfrac{\log(Aq_{68}/\epsilon(\mu))}{\log B}< 32 <n-m.$$
So

$$n-m\leq 32.$$
Substituting this upper bound for $n-m$ in the inequality \eqref{eq:8}, we obtain $k< 1.32\cdot 10^{16}.$ Now let 

$$z_2:=n\log 10-k\log \alpha+\log((a_1-a_210^{m-n})/9\cdot 2\sqrt{2}).$$

From the inequality \eqref{eq:7}, we have that 

$$|x|=|e^{z_1}-1|<\frac{4}{\alpha^k}<\frac{1}{10}$$
for $k\geq 150.$ Choosing $a=0.1$, we get from Lemma \ref{weber} the inequality

$$|z_2|=|\log(x+1)|<\frac{\log(10/9)}{1/10}\cdot \frac{4}{\alpha^k}<4.22\cdot \alpha^{-k}.$$
Thus it follows that 

$$0<\Big|n\log 10-k\log \alpha+\log((a_1-a_210^{m-n})/9\cdot 2\sqrt{2})\Big|<4.22\cdot \alpha^{-k}$$
Dividing both sides by $\log\alpha$, we obtain 

\begin{equation}
\label{eq:10}
0<0<\Big|n\left(\frac{\log 10}{\log\alpha}\right)-k+\left(\frac{\log((a_1-a_210^{m-n})/9\cdot 2\sqrt{2})}{\log\alpha}\right)\Big|<4.79\cdot \alpha^{-k}.
\end{equation}

Put $\gamma:=\frac{\log 10}{\log \alpha}\not\in \mathbb{Q}$ and $M= 1.32\cdot 10^{16}.$ Then we found that $q_{42}$, the denominator of the $42$-th convergent of $\gamma$ is greater than $6M$. Now take 

$$\mu:=\frac{\log((a_1-a_210^{m-n})/9\cdot 2\sqrt{2})}{\log\alpha}.$$
In this case, considering the fact that $1\leq a_1,a_2\leq 9$ and $2\leq n-m\leq 32$, a quick computation with Sagemath gives us that 

$$\epsilon(\mu):=||\mu q_{44}||-M||\gamma q_{42}|| > 0.0002269.$$
Let $A=4.8$, $B=10$ and $w=k$ in Lemma \ref{Dujjella}. Thus, using Sagemath, we can say that the inequality \eqref{eq:10} has no solution if 

$$ \dfrac{\log(Aq_{42}/\epsilon(\mu))}{\log B}<23<k.$$
This implies that $k\leq 23$. This contradicts our assumption that $k\geq 150$. This completes the proof.

\section*{Proof of Theorem \ref{Main1}}

The proof is similar to that of Theorem \ref{Main}. We may sometimes omit some details.
\medskip

Assume that the equation  \eqref{Main1} holds. Let $1\leq k \leq 149$ and $n\geq 2$. Then, by using Sagemath, we obtain only the solutions listed in Theorem \ref{Main1b}.
\medskip

From now, we assume that $k\geq 150.$ If $n=m$, then it follows that $a_1>a_2$, which means that $Q_k$ is a repdigit. But the largest repdigit in $Q_k$ is $6$. Thus we get a contradiction since $k\geq 150.$ 
\medskip

The case $n-m=1$  has been completed solved in \cite{ALT} by Adedji at al.  They founded that the largest Pell-Lucas number that can be expressed as a  concatenations of two repdigits is $Q_5=82.$

Assume that the equation \eqref{Main1b} have solutions and let $k\geq 150$, $n-m\geq 2$. Then, using \eqref{eq:siezQn}, we get the inequality

$$\frac{\alpha^{2n}}{20}<\frac{10^{n-1}}{2}<10^{n-1}-10^{m}<a_1\frac{10^n-1}{9}-a_2\frac{10^m-1}{9}=Q_k \leq \alpha^{k+1}, $$
giving that $2n<k+5.$ In particular $n<k+2.$ On the other hand, we rewrite equation \eqref{Main1b} as

$${\alpha^k + \beta^k}= a_1\frac{10^n-1}{9}-a_2\frac{10^m-1}{9}$$
to obtain

\begin{equation}
\label{eq:1b}
{9\alpha^k}-a_110^n={-9\beta^k}-a_210^m-(a_1-a_2).
\end{equation}
Taking absolute value of both sides of equation \eqref{eq:1b}, we obtain 
\begin{equation}
\label{eq:2b}
\Big|{9\alpha^k}-a_110^n\Big|\leq {9|\beta|^k}+a_210^m+|a_1-a_2|.
\end{equation}
Dividing both sides of \eqref{eq:2b} by $a_110^n$, we obtain

\begin{align*}
\Big|{9\cdot 10^{-n}a_1^{-1}\alpha^k}-1\Big|& \leq \frac{9|\beta|^k}{a_110^n}+\frac{a_210^m}{a_110^n}+\frac{|a_1-a_2|}{a_110^n}\\
&\leq \frac{9|\beta|^k}{10^{n-m+1}}+\frac{9}{10^{n-m}}+\frac{8}{10^{n-m+1}}.
\end{align*}
This implies that 

\begin{equation}
\label{eq:3b}
\Big|\frac{9\cdot 10^{-n}\alpha^k}{a_1}-1\Big| \leq \frac{21}{10^{n-m}}.
\end{equation}

Now we apply Theorem \ref{Matveev11} to the left-hand side of the above inequality with $(\gamma_1,\gamma_2,\gamma_3)=(\alpha,10,9/a_1
)$ and $(b_1,b_2,b_3)=(k,-n,1)$. Note that $\gamma_1,\gamma_2$ and $\gamma_3$ are positive real numbers and elements of the field $\mathbb{K}=\mathbb{Q}(\sqrt{2})$. Therefore the degree of the field $\mathbb{K}$ is equal to $D=2$.
Put 

$$\Lambda_1=\frac{9\cdot 10^{-n}\alpha^k}{a_1}-1.$$
If $\Lambda_1=0$, then $\alpha^k=10^na_1/9 \in \mathbb{Q}$, which is a contradiction.  Therefore $\Lambda_1\neq 0$. Using the properties of the logarithmic height in \eqref{eq:height}, we have that

$$h(\gamma_1)=h(\alpha)=\frac{\log\alpha}{2}\leq\frac{0.882}{2}, \quad h(\gamma_2)=\log10$$
and 

$$h(\gamma_3)=h(9/a_1)< h(9)+h(a_1) <2\log{9}.$$ We can take $A_1=\log(\alpha)=0.882, A_2=4.7$ and $A_3=9$. Since $n<k$ and $B\geq \max(n,k,1)$ we can take $B=k+2$. Thus, we obtain by inequality \eqref{eq:3b} and  Theorem \ref{Matveev11} that 

$$21\cdot 10^{m-n}>\Lambda_1>exp(C\cdot (1+\log2)\cdot (1+\log(k+2))(0.882)(4.7)(9))$$
where $C=-1.4\cdot 30^6\cdot 3^{4.5}\cdot 2^2.$ It follows that 

\begin{equation}
\label{eq:4b}
m-n< 1.8\cdot 10^{13}\cdot (1+\log(k+2)).
\end{equation}

Rearranging equation \eqref{Main1b} 	as 

\begin{equation}
\label{eq:5b}
{\alpha^k}-\frac{a_110^n-a_210^m}{9}={\beta^k}-\frac{(a_1-a_2)}{9}
\end{equation}
and taking absolute value of both sides of \eqref{eq:5b}, we get

\begin{equation}
\label{eq:6b}
\Big|{\alpha^k}-\frac{a_110^n-a_210^m}{9}\Big|\leq {\Big|\beta\Big|^k}+\frac{\Big|a_1-a_2\Big|}{9}.
\end{equation}
Dividing both sides of the above inequality by ${\alpha^k}$ we obtain 

\begin{equation}
\label{eq:7b}
\Big|1-{(a_1-a_210^{m-n})\cdot 10^n\cdot \alpha^{-k}}\cdot 9^{-1}\Big|\leq \frac{1}{\alpha^{2k}}+\frac{16}{9\alpha^k}<\frac{4}{\alpha^k}.
\end{equation}
Now, we can apply again Theorem \ref{Matveev11} to the above inequality  with

$$(\gamma_1,\gamma_2,\gamma_3)=(\alpha,10,(a_1-a_210^{m-n})/9)$$ and 

$$(b_1,b_2,b_3)=(-k,n,1).$$

Note that $\gamma_1,\gamma_2$ and $\gamma_3$ are positive real numbers and elements of the field $\mathbb{K}=\mathbb{Q}(\sqrt{2})$. Therefore the degree of the field $\mathbb{K}$ is equal to $D=2$. Let 

$$\Lambda_2= 1-\frac{(a_1-a_210^{m-n})\cdot 10^n\cdot \alpha^{-k}}{9}.$$
If $\Lambda=0$, then get that $\alpha^{2k}\in \mathbb{Q}$, which is false for $k>0$. Using the properties of the logarithmic height in \eqref{eq:height}, we have that 

$$h(\gamma_1)=h(\alpha)=\frac{\log\alpha}{2}\leq\frac{0.882}{2}, \quad h(\gamma_2)=\log10$$
and 

\begin{align*}
h(\gamma_3)&=h((a_1-a_210^{m-n})/9)\\
&\leq h(9)+h(a_1)+h(a_2)+(n-m)h(10)+\log2\\
&\leq 8.33 +(n-m)\log10.
\end{align*}
By Equation \eqref{eq:4b}, we get:
$$h(\gamma_3) \leq 35\cdot 10^{13}\cdot (1+\log(k+2)).$$
So we can take $A_1=\log(\alpha)=0.9, A_2=4.7$ and $A_3=35\cdot 10^{13}\cdot (1+\log(k+2))$. Since $n<k$ and $B\geq \max(n,k,1)$, we can take $B=k$. Thus, we obtain by inequality \eqref{eq:7b} and  Theorem \ref{Matveev11} that 

$$4\cdot \alpha^{-k}>\Lambda_2>\exp(C\cdot (1+\log2)\cdot (1+\log(k+2))(0.9)(4.6)(35\cdot 10^{13}\cdot (1+\log(k+2)))),$$
where $C=-1.4\cdot 30^6\cdot 3^{4.5}\cdot 2^2.$ By a simple computation, it follows that 

\begin{equation}
\label{eq:8b}
k < 1.7\cdot 10^{27}\cdot (1+\log(k+2))^2.
\end{equation}
Now using the nequalities  \eqref{eq:8b} and \eqref{eq:4b}, a computer search with SageMath gives us that $k<9\cdot 10^{30}$ and 
$n-m<13\cdot 10^{14}$.

\subsection{Reducing the Bounding of $ k $}
Now let us reduce the upper bound on $k$. Firstly, we determine  a suitable upper bound on $n-m$, and later we use Lemma \ref{Dujjella} to conclude that $n$ must be  smaller than $150$. Let 

$$z_1:= k\log \alpha-n\log 10 +\log(9/a_1).$$
From the inequality \eqref{eq:3b}, we have that 
$$|x|=|e^{z_1}-1|\leq \frac{21}{10^{n-m}}<\frac{21}{10^2}$$
for $n-m\geq 2.$ Choosing $a:=0.1$, we get the inequality

$$z_1=|\log(x+1)|<\frac{-\log(1-21/10^2)}{21/10^2}\cdot\frac{21}{10^{n-m}}<21\cdot 10^{m-n}$$
by Lemma \ref{weber}. Thus it follows that 

$$0<|k\log \alpha-n\log 10 +\log(9/a_1)|<21\cdot 10^{m-n}.$$
\medskip

Dividing this inequality by $\log 10$, we obtain

\begin{equation}
\label{eq:9b}
\Big|k\left(\frac{\log \alpha}{\log 10}\right)-n +\left(\frac{\log(9/a_1)}{\log 10}\right)\Big|<10^{m-n+1}.
\end{equation}
We can take $\gamma:=\frac{\log \alpha}{\log 10}\not\in \mathbb{Q}$ and $M= 13\cdot 10^{14}.$ Then we found that $q_{70}$, the denominator of the $70$-th convergent of $\gamma$ is greater than $6M$. Now take 

$$\mu:=\frac{\log(9/a_1)}{\log 10}.$$
In this case, considering the fact that $1\leq a_1< 9$, a quick computation with Sagemath gives us that

$$\epsilon(\mu):=||\mu q_{38}||-M||\gamma q_{38}||>0.02249.$$
Let $A=10$, $B=10$ and $w=n-m$ in Lemma \ref{Dujjella}. Thus, using Sagemath, we can say that the inequality \eqref{eq:9b} has no solution if  

$$ \dfrac{\log(Aq_{38}/\epsilon(\mu))}{\log B}< 19 <n-m.$$
So

$$n-m\leq 19.$$

When $a_1 = 9$,  the parameter $\mu = 0$; and the resulting inequality from \eqref{eq:9b} has the shape
\begin{center}
	$0 < \vert a \gamma - b \vert < 10\cdot 10^{m-n}$,
\end{center}
with $\gamma$ being an irrational number and $a,b \in \mathbb{Z}$. So, one can appeal to the known properties of the convergents of the continued fractions to obtain a nontrivial lower bound for

\begin{center}
	$\vert a \gamma - b \vert$.
\end{center}

Let us see.
When $a_1= 9$, from \eqref{eq:9b}, we get that

\begin{equation}
	\label{eq:9b1}
	0  <  k\gamma - n  < \frac{10}{10^{n-m}}.
\end{equation}
Let $[a_0 ,a_1 ,a_2 ,a_3 ,a_4,a_5,a_6,a_7,...] = [0; 2, 1, 1, 2, 2\ldots]$ be the continued fraction expansion of $\gamma$, and let
denote $p_n /q_n$ its $n$th convergent. Recall also that $k < 9\cdot10^{30}$.

Furthermore, $a_N := \max\lbrace{a_i: i = 0,1,\ldots,38}\rbrace = a_{24} = 39$. So, from the known properties
of continued fractions, we obtain that
\begin{equation}
	\label{eq:9b2}
	\vert k\gamma - n\vert  >  \frac{1}{(a_N + 2)k}.
\end{equation}
Comparing estimates \eqref{eq:9b1} and \eqref{eq:9b2}, we get right away that
\begin{equation}
	\label{eq:9b3} 
	10^{n-m}  <  10\cdot 41\cdot k < 3.69\cdot 10^{33},
\end{equation}
leading to $n-m < 34$. Consequently,
$n-m< 34$ always holds.

 Now let 

$$z_2:=n\log 10-k\log \alpha+\log((a_1-a_210^{m-n})/9).$$

From the inequality \eqref{eq:7b}, we have that 

$$|x|=|e^{z_1}-1|<\frac{4}{\alpha^k}<\frac{1}{10}$$
for $k\geq 150.$ Choosing $a=0.1$, we get from Lemma \ref{weber} the inequality

$$|z_2|=|\log(x+1)|<\frac{\log(10/9)}{(1/10)}\cdot \frac{4}{\alpha^k}<5\cdot \alpha^{-k}.$$
Thus it follows that 

$$0<\Big|n\log 10-k\log \alpha+\log((a_1-a_210^{m-n})/9)\Big|<5\cdot \alpha^{-k}$$

Dividing both sides by $\log\alpha$, we obtain 

\begin{equation}
\label{eq:10b}
0<0<\Big|n\left(\frac{\log 10}{\log\alpha}\right)-k+\left(\frac{\log((a_1-a_210^{m-n})/9)}{\log\alpha}\right)\Big|<5.7\cdot \alpha^{-k}
\end{equation}

Put $\gamma:=\frac{\log 10}{\log \alpha}\not\in \mathbb{Q}$ and $M= 9\cdot 10^{30}.$ Then we found that $q_{68}$, the denominator of the $68$-th convergent of $\gamma$ is greater than $6M$. Now take 

$$\mu:=\frac{\log((a_1-a_210^{m-n})/9)}{\log\alpha}.$$
In this case, considering the fact that $1\leq a_1,a_2\leq 9$ and $2\leq n-m\leq 34$, a quick computation with SageMath gives us that 

$$\epsilon(\mu):=||\mu q_{69}||-M||\gamma q_{69}||> 0.0007213.$$
Let $A=5.7$, $B=10$ and $w=k$ in Lemma \ref{Dujjella}. Thus, using Sagemath, we can say that the inequality \eqref{eq:10b} has no solution if 

$$ \dfrac{\log(Aq_{69}/\epsilon(\mu))}{\log B}<36<k.$$
This implies that $k\leq 36.$ This contradicts our assumption that $k>150$. This completes the proof.

\medskip

\textbf{Author Contributions} The authors contributed equally in writing the final version of this article. All authors
read and approved the final manuscript.

\section*{Declarations}
\textbf{Conflict of Interest:} The authors declare that they have no conflict of interest.

\section*{Acknowledgements}

We are very grateful to the referee for the careful reading of the paper and for his comments and detailed suggestions which helped us to improve considerably the manuscript. 

 \end{document}